# Sum of cubes: Old proofs suggest new $q-$ analogues

Johann Cigler

Fakultät für Mathematik, Universität Wien

johann.cigler@univie.ac.at

**Abstract**

We prove a new $q-$ analogue of Nicomachus's theorem about the sum of cubes and some related results.

**1. Introduction**

In [1], [5] and [8] some $q-$ analogues of the well-known formula

$$1^3 + 2^3 + \cdots + n^3 = \binom{n+1}{2}^2 \tag{1}$$

have been found. Some information about this formula which is sometimes called Nicomachus's theorem is given in [7].

In this note I propose another $q-$ analogue which is inspired by an old result of C. Wheatstone [6].
He observed that the odd numbers can be grouped in such a way that the identities

$$\begin{aligned} 1 &= 1^3 \\ 3+5 &= 2^3 \\ 7+9+11 &= 3^3 \\ 13+15+17+19 &= 4^3 \\ &\cdots \end{aligned} \tag{2}$$

hold, which implies that

$$1^3 + 2^3 + \cdots + n^3 = 1 + 3 + \cdots + \left(2\binom{n+1}{2} - 1\right).$$

(1) follows from the well-known formula

$$1 + 3 + \cdots + (2n-1) = n^2. \tag{3}$$

Identity (1) is usually stated in the form

$$1^3 + 2^3 + \cdots + n^3 = (1 + 2 + \cdots + n)^2. \qquad (4)$$

This is an immediate consequence of

$$1 + 2 + \cdots + n = \binom{n+1}{2}. \qquad (5)$$

A simple "proof without words" of (4) has been given in [4], which I reproduce here:

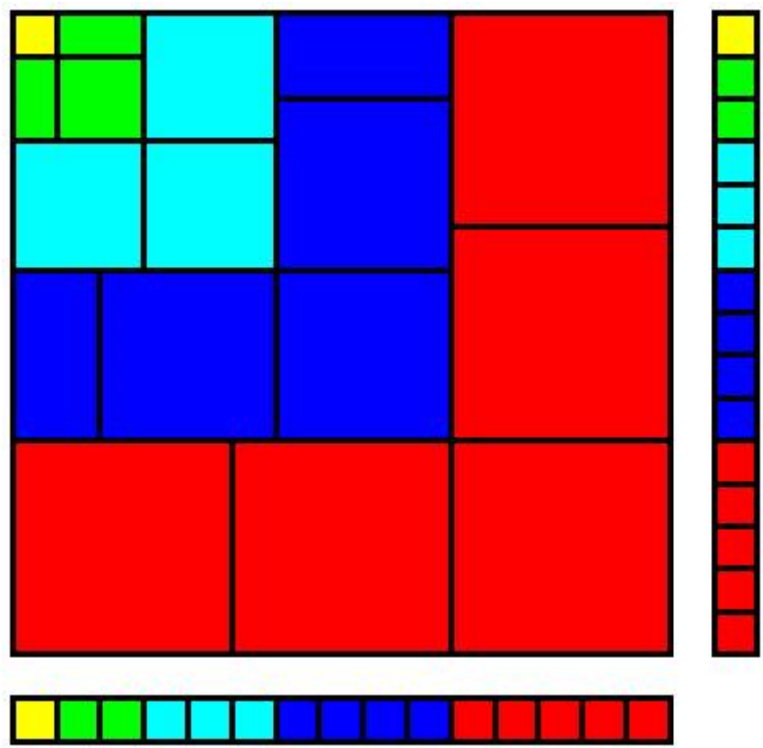

The simplest computational proof of (1) uses the trivial identity $\binom{n+1}{2}^2 - \binom{n}{2}^2 = n^3$ which gives the telescoping sum

$$1^3 + 2^3 + \cdots + n^3 = \binom{2}{2}^2 - \binom{1}{2}^2 + \binom{3}{2}^2 - \binom{2}{2}^2 + \cdots + \binom{n+1}{2}^2 - \binom{n}{2}^2 = \binom{n+1}{2}^2 - \binom{1}{2}^2 = \binom{n+1}{2}^2.$$

## 2. $q$ − analogues.

As usual we let $[n]_q = \dfrac{1-q^n}{1-q} = 1+q+\cdots+q^{n-1}$ and $\begin{bmatrix} n \\ k \end{bmatrix}_q = \dfrac{[n]_q[n-1]_q\cdots[n-k+1]_q}{[1]_q[2]_q\cdots[k]_q}$ for $0 < q < 1$. It is clear that $\lim\limits_{q \to 1}[n]_q = n$ and $\lim\limits_{q \to 1}\begin{bmatrix} n \\ k \end{bmatrix}_q = \binom{n}{k}$.

K.C. Garrett and K. Hummel [1] derived a combinatorial proof of the $q$ − analogue

$$\sum_{k=1}^{n} q^{k-1}[k]_q^2 \frac{[k-1]_q + [k+1]_q}{[2]_q} = \begin{bmatrix} n+1 \\ 2 \end{bmatrix}_q^2, \tag{6}$$

S.O. Warnaar [5] proposed the identity

$$\sum_{k=1}^{n} q^{2n-2k}[k]_q^2[k]_{q^2} = \begin{bmatrix} n+1 \\ 2 \end{bmatrix}_q^2 \tag{7}$$

and G. Zhao and H. Feng [8] gave a combinatorial interpretation of

$$\sum_{k=0}^{n} q^{4(n-k)}[k]_q^2 \frac{1+q^2-2q^{k+1}}{1-q^2} = \begin{bmatrix} n+1 \\ 2 \end{bmatrix}_q^2. \tag{8}$$

For a proof by induction observe that

$$\begin{bmatrix} n+1 \\ 2 \end{bmatrix}_q^2 - q^{2a}\begin{bmatrix} n \\ 2 \end{bmatrix}_q^2 = \frac{(1-q^n)^2}{(1-q)^2(1-q^2)^2}\left((1-q^{n+1})^2 - (q^a - q^{n-1+a})^2\right)$$

$$= \frac{(1-q^n)^2}{(1-q)^2(1-q^2)^2}(1-q^{n+1} - (q^a - q^{n-1+a}))(1-q^{n+1} + (q^a - q^{n-1+a}))$$

implies

$$\begin{bmatrix} n+1 \\ 2 \end{bmatrix}_q^2 - \begin{bmatrix} n \\ 2 \end{bmatrix}_q^2 = q^{n-1}[n]_q^2 \frac{[n-1]_q + [n+1]_q}{[2]_q},$$

$$\begin{bmatrix} n+1 \\ 2 \end{bmatrix}_q^2 - q^2\begin{bmatrix} n \\ 2 \end{bmatrix}_q^2 = [n]_q^2[n]_{q^2}$$

and

$$\begin{bmatrix} n+1 \\ 2 \end{bmatrix}_q^2 - q^4 \begin{bmatrix} n \\ 2 \end{bmatrix}_q^2 = [n]_q^2 \frac{1+q^2-2q^{n+1}}{1-q^2}.$$

A similar identity is

$$\left[\binom{n+1}{2}\right]_q^2 - q^n \left[\binom{n}{2}\right]_q^2 = [n]_q [n^2]_q, \qquad (9)$$

which follows from

$$\left(1-q^{\binom{n+1}{2}}\right)^2 - q^n\left(1-q^{\binom{n}{2}}\right)^2 = 1-2q^{\binom{n+1}{2}}+q^{2\binom{n+1}{2}}-q^n\left(1-2q^{\binom{n}{2}}+q^{2\binom{n}{2}}\right) = (1-q^n)(1-q^{n^2}).$$

This identity gives the telescoping sum

$$\sum_{j=1}^{n} q^{\binom{n+1}{2}-\binom{j+1}{2}} [j]_q [j^2]_q = q^{\binom{n+1}{2}} \sum_{j=1}^{n} q^{-\binom{j+1}{2}} \left( \left[\binom{j+1}{2}\right]_q^2 - q^j \left[\binom{j}{2}\right]_q^2 \right)$$

$$= q^{\binom{n+1}{2}} \left( q^{-\binom{2}{2}} \left[\binom{2}{2}\right]_q^2 - q^{-\binom{1}{2}} \left[\binom{1}{2}\right]_q^2 + q^{-\binom{3}{2}} \left[\binom{3}{2}\right]_q^2 - q^{-\binom{2}{2}} \left[\binom{2}{2}\right]_q^2 + \cdots + q^{-\binom{n+1}{2}} \left[\binom{n+1}{2}\right]_q^2 - q^{-\binom{n}{2}} \left[\binom{n}{2}\right]_q^2 \right)$$

$$= \left[\binom{n+1}{2}\right]_q^2 - q^{\binom{n+1}{2}} q^{-\binom{1}{2}} \left[\binom{1}{2}\right]_q^2 = \left[\binom{n+1}{2}\right]_q^2.$$

Thus we have obtained our main result

**Theorem 1**

$$\sum_{j=1}^{n} q^{\binom{n+1}{2}-\binom{j+1}{2}} [j]_q [j^2]_q = \left[\binom{n+1}{2}\right]_q^2. \qquad (10)$$

We now give two further proofs which generalize the beautiful proofs which we have sketched in the introduction and which originally led to this $q-$ analogue.

We start with the following well-known (cf. [3]) $q-$ analogue of (3)

$$\sum_{k=1}^{n} q^{n-k}[2k-1]_q = [n]_q^2. \tag{11}$$

A computational proof uses the fact that

$$[n]_q^2 - q[n-1]_q^2 = \frac{q^{2n} - 2q^n + 1 - q^{2n-1} + 2q^n - q}{(q-1)^2} = \frac{q^{2n-1} - 1}{q-1} = [2n-1]_q.$$

But formula (11) has also a nice combinatorial interpretation.

Consider the squares $S_n = \{(i,j)\}_{0 \le i,j < n}$ and associate with each point $(i,j) \in S_n$ the weight $w(i,j) = q^i q^{n-1-j}$. The weight of $S_n$ is

$$w(S_n) = \sum_{i=0}^{n-1} q^i \sum_{j=0}^{n-1} q^{n-1-j} = [n]_q^2. \tag{12}$$

The square $S_n$ is the union of the hooks
$$h_k = \{(k-1,0),(k-1,1),\cdots,(k-1,k-1),(k-2,k-1),\cdots,(0,k-1)\}, \ 1 \le k \le n.$$

The weight of the hook $h_k$ is
$$w(h_k) = q^n\left(q^{k-2} + q^{k-3} + \cdots + q^0 + q^{-1} + \cdots + q^{-k}\right) = q^{n-k}\left(1 + q + \cdots + q^{2k-2}\right) = q^{n-k}[2k-1]_q.$$

The point $(n,0)$ will be called the base-point of the square.

As an example consider the case $n=6$ in matrix notation.

$$\begin{pmatrix} q^5 & q^6 & q^7 & q^8 & q^9 & q^{10} \\ q^4 & q^5 & q^6 & q^7 & q^8 & q^9 \\ q^3 & q^4 & q^5 & q^6 & q^7 & q^8 \\ q^2 & q^3 & q^4 & q^5 & q^6 & q^7 \\ q & q^2 & q^3 & q^4 & q^5 & q^6 \\ 1 & q & q^2 & q^3 & q^4 & q^5 \end{pmatrix} \tag{13}$$

Here we have $w(h_1) = q^5$, $w(h_2) = q^4[3]_q$, $w(h_2) = q^3[5]_q, \cdots$.

## Second proof of Theorem 1

We first observe the nice $q-$analogues

$$1 = [1]_q[1^2]_q$$
$$q[3]_q + [5]_q = [2]_q[2^2]_q$$
$$q^2[7]_q + q[9]_q + [11]_q = [3]_q[3^2]_q$$

and more generally

$$\sum_{j=1}^{n} q^{n-j}\left[n^2 - n + 2j - 1\right]_q = [n]_q[n^2]_q \qquad (14)$$

of formulae (2).

They follow from the identity

$$\sum_{j=1}^{n} q^{n-j}\left[n^2 - n + 2j - 1\right]_q = \frac{1}{1-q}\sum_{j=1}^{n} q^{n-j}\left(1 - q^{n^2-n+2j-1}\right) = \frac{1}{1-q}\left(\sum_{j=0}^{n-1} q^{n-1-j} - \sum_{j=0}^{n-1} q^{n^2+j}\right)$$
$$= \frac{[n]_q - q^{n^2}[n]_q}{1-q} = [n]_q[n^2]_q.$$

Using (11) we get the desired result

$$\sum_{j=1}^{n} q^{\binom{n+1}{2} - \binom{j+1}{2}} [j]_q[j^2]_q = \sum_{j=1}^{n} q^{\binom{n+1}{2} - \binom{j+1}{2}} \sum_{k=1}^{j} q^{j-k}\left[j^2 - j + 2k - 1\right]_q = \sum_{j=1}^{n} q^{\binom{n+1}{2}} \sum_{k=\binom{j}{2}+1}^{\binom{j+1}{2}} q^{-k}\left[2k - 1\right]_q$$

$$= \sum_{k=1}^{\binom{n+1}{2}} q^{\binom{n+1}{2} - k} [2k - 1]_q = \left[\binom{n+1}{2}\right]_q^2.$$

## Third proof of Theorem 1

A combinatorial proof can also be given along the lines of the above "proof without words".

For odd $j$ the union $R_j$ of the $j$ hooks $h_{\binom{j}{2}+1}, \cdots, h_{\binom{j+1}{2}}$ is the union of $j$ squares of side-length $j$ whose base points have weight $q^{\binom{n+1}{2} - \binom{j+1}{2}} q^{mj}$, $0 \leq m < j$.

The weight of these squares is

$$q^{\binom{n+1}{2}-\binom{j+1}{2}}[j]_q^2\left(1+q^j+\cdots q^{(j-1)j}\right)=q^{\binom{n+1}{2}-\binom{j+1}{2}}[j]_q^2[j]_{q^j}=q^{\binom{n+1}{2}-\binom{j+1}{2}}[j]_q\frac{1-q^j}{1-q}\frac{1-q^{j^2}}{1-q^j}$$

$$=q^{\binom{n+1}{2}-\binom{j+1}{2}}[j]_q[j^2]_q.$$

For $j=2\ell$ the union $R_j$ of the hooks $h_{\binom{j}{2}+1},\ldots,h_{\binom{j+1}{2}}$ is the union of $j-1$ squares whose base-points have weights $q^{\binom{n+1}{2}-\binom{j+1}{2}}q^{\ell+2\ell m}$, $0\le m\le j-2$, and of two rectangles with side lengths $j$ and $\ell$ as in the blue region in the above figure.

The weight of the uppermost rectangle is $q^{\binom{n+1}{2}+\binom{j}{2}-\ell}[j]_q[\ell]_q$ and the weight of the leftmost rectangle is $q^{\binom{n+1}{2}-\binom{j+1}{2}}[j]_q[\ell]_q$.

Thus the total weight of this region is

$$q^{\binom{n+1}{2}-\binom{2\ell+1}{2}}\sum_{m=0}^{2\ell-2}q^{\ell+2\ell m}[2\ell]_q^2+q^{\binom{n+1}{2}+\binom{2\ell}{2}-\ell}[2\ell]_q[\ell]_q+q^{\binom{n+1}{2}-\binom{2\ell+1}{2}}[2\ell]_q[\ell]_q$$

$$=q^{\binom{n+1}{2}-\binom{2\ell+1}{2}}[2\ell]_q[\ell]_q\left(\sum_{m=0}^{2\ell-2}q^{\ell+2\ell m}\frac{1-q^{2\ell}}{1-q^\ell}+q^{4\ell^2-\ell}+1\right)$$

$$=q^{\binom{n+1}{2}-\binom{2\ell+1}{2}}[2\ell]_q[\ell]_q\left(\frac{q^\ell-q^{4\ell^2-\ell}}{1-q^\ell}+q^{4\ell^2-\ell}+1\right)=q^{\binom{n+1}{2}-\binom{2\ell+1}{2}}[2\ell]_q\frac{1-q^\ell}{1-q}\frac{\left(1-q^{4\ell^2}\right)}{1-q^\ell}$$

$$=q^{\binom{n+1}{2}-\binom{2\ell+1}{2}}[2\ell]_q[4\ell^2]_q.$$

Thus we get $\sum_{j=1}^n w(R_j)=w\left(S_{\binom{n+1}{2}}\right)=\left[\binom{n+1}{2}\right]_q^2$ and thus again (10).

## 3. Related results

Let us note some related results. For any sequence of positive integers $a(n)$ the sums

$$\sum_{j=1}^n q^{\sum_{i=1}^{j-1}a(i)}[a(j)]_q=\left[\sum_{i=1}^n a(i)\right]_q \tag{15}$$

and

$$\sum_{j=1}^{n} q^{\sum_{i=j+1}^{n} a(i)} [a(j)]_q = \left[\sum_{i=1}^{n} a(i)\right]_q \qquad (16)$$

are $q-$ analogues of $\sum_{i=1}^{n} a(i)$.

The proofs are obvious because

$$\left[\sum_{i=1}^{n} a(i)\right]_q + q^{\sum_{i=1}^{n} a(i)} [a(n+1)]_q = \frac{1-q^{\sum_{i=1}^{n} a(i)}}{1-q} + \frac{q^{\sum_{i=1}^{n} a(i)}\left(1-q^{a(n+1)}\right)}{1-q} = \frac{1-q^{\sum_{i=1}^{n+1} a(i)}}{1-q} = \left[\sum_{i=1}^{n+1} a(i)\right]_q$$

and

$$q^{a(n+1)} \left[\sum_{i=1}^{n} a(i)\right]_q + [a(n+1)]_q = \frac{q^{a(n+1)} - q^{\sum_{i=1}^{n+1} a(i)}}{1-q} + \frac{1-q^{a(n+1)}}{1-q} = \frac{1-q^{\sum_{i=1}^{n+1} a(i)}}{1-q} = \left[\sum_{i=1}^{n+1} a(i)\right]_q.$$

By choosing $a(n) = n^3$ we get the following $q-$ analogues of (1):

**Theorem 2**

$$\sum_{j=1}^{n} q^{\binom{j}{2}^2} [j^3]_q = \left[\binom{n+1}{2}^2\right]_q \qquad (17)$$

*and*

$$\sum_{j=1}^{n} q^{\sum_{i=j+1}^{n} i^3} [j^3]_q = \left[\binom{n+1}{2}^2\right]_q. \qquad (18)$$

From the recurrence relations for the $q-$ binomial coefficients

$$\begin{bmatrix} n+1 \\ k \end{bmatrix}_q = q^k \begin{bmatrix} n \\ k \end{bmatrix}_q + \begin{bmatrix} n \\ k-1 \end{bmatrix}_q$$

$$\begin{bmatrix} n+1 \\ k \end{bmatrix}_q = \begin{bmatrix} n \\ k \end{bmatrix}_q + q^{n-k+1} \begin{bmatrix} n \\ k-1 \end{bmatrix}_q$$

we get the well-known formulae

$$\sum_{j=1}^{n} q^{j-1} \begin{bmatrix} j \\ k \end{bmatrix}_q = q^{k-1} \begin{bmatrix} n+1 \\ k+1 \end{bmatrix}_q \tag{19}$$

and

$$\sum_{j=1}^{n} q^{(k+1)(n-j)} \begin{bmatrix} j \\ k \end{bmatrix}_q = \begin{bmatrix} n+1 \\ k+1 \end{bmatrix}_q. \tag{20}$$

For $q=1$ these sums can be used to compute $\sum_{k=1}^{n} j^k$. For example from $n^3 = 6\binom{n}{3} + 6\binom{n}{2} + n$

we get $\sum_{j=1}^{n} j^3 = 6\binom{n+1}{4} + 6\binom{n+1}{3} + \binom{n+1}{2} = \frac{n^2(n+1)^2}{4} = \binom{n+1}{2}^2$.

Unfortunately in general the sums $\sum_{k=1}^{n} q^{j-1}[j]^k$ don't have a simple expression.

For example from $[n]_q^3 = q^3 [2]_q [3]_q \begin{bmatrix} n \\ 3 \end{bmatrix}_q + [2]_q [3]_q \begin{bmatrix} n \\ 2 \end{bmatrix}_q + q^{n-1}[n]_q$ we get

$$\sum_{j=1}^{n} q^{j-1}[j]_q^3 = q^5 [2]_q [3]_q \begin{bmatrix} n+1 \\ 4 \end{bmatrix}_q + q[2]_q [3]_q \begin{bmatrix} n+1 \\ 3 \end{bmatrix}_q + \sum_{j=1}^{n} q^{2(j-1)}[j]_q.$$

Now

$$\sum_{j=1}^{n} q^{2(j-1)}[j]_q = \sum_{j=1}^{n} q^{(j-1)}[j]_q + \sum_{j=1}^{n} q^{(j-1)}\left(q^{(j-1)} - 1\right)[j]_q = \begin{bmatrix} n+1 \\ 2 \end{bmatrix}_q - (1-q^2)\sum_{j=1}^{n} q^{(j-1)} \begin{bmatrix} j \\ 2 \end{bmatrix}_q$$

$$= \begin{bmatrix} n+1 \\ 2 \end{bmatrix}_q - q(1-q^2)\begin{bmatrix} n+1 \\ 3 \end{bmatrix}_q = \frac{[n+1]_q [n]_q}{[2][3]}\left(1+q+q^2 - q(1+q)(1-q^{n-1})\right) = \begin{bmatrix} n+1 \\ 2 \end{bmatrix}_q \frac{1+q^n+q^{n+1}}{1+q+q^2}.$$

Thus we get

$$\sum_{j=1}^{n} q^{j-1}[j]_q^3 = q^5[2]_q[3]_q \begin{bmatrix} n+1 \\ 4 \end{bmatrix}_q + q[2]_q[3]_q \begin{bmatrix} n+1 \\ 3 \end{bmatrix}_q + \begin{bmatrix} n+1 \\ 2 \end{bmatrix}_q \frac{1+q^n+q^{n+1}}{1+q+q^2}. \quad (21)$$

This formula cannot be simplified.

A curious generalization of (2) is due to P. Luthy [2].

He observed that

$$\begin{aligned} 1 &= 1^5 \\ 5+7+9+11 &= 2^5 \\ 19+21+23+25+27+29+31+33+35 &= 3^5 \\ &\ldots \end{aligned} \quad (22)$$

and more generally

$$\sum_{j=1}^{n^k} \left( n^{k+1} - n^k + 2j - 1 \right) = n^{2k+1}. \quad (23)$$

To see this observe that $\dfrac{\left(n^{k+1} - n^k + 2j - 1\right) + \left(n^{k+1} - n^k + 2(n^k - j) + 1\right)}{2} = n^{k+1}$ for $1 \le j \le n^k$.

Similar results also hold for the $q-$ analogues of (23).

$$\sum_{j=1}^{n^k} q^{n^k - j} \left[ n^{k+1} - n^k + 2j - 1 \right]_q = \left[ n^k \right]_q \left[ n^{k+1} \right]_q. \quad (24)$$

For

$$\sum_{j=1}^{n^k} q^{n^k - j} \left[ n^{k+1} - n^k + 2j - 1 \right]_q = \frac{1}{1-q} \sum_{j=1}^{n^k} q^{n^k - j} \left( 1 - q^{n^{k+1} - n^k + 2j - 1} \right) = \frac{1}{1-q} \sum_{j=0}^{n^k - 1} q^j - \frac{q^{n^{k+1}}}{1-q} \sum_{j=0}^{n^k - 1} q^j$$

$$= \frac{1}{(1-q)^2} \left( 1 - q^{n^k} \right)\left( 1 - q^{n^{k+1}} \right) = \left[ n^k \right]_q \left[ n^{k+1} \right]_q.$$

## 4. More $q-$ analogues of identities by C. Wheatstone [6].

### 4.1

A $q-$ analogue of $\sum_{j=0}^{n-1}(n+1+2j)=2n^2$ is

$$\sum_{j=0}^{n-1}[n+1+2j]_q\, q^{n-1-j}=[2]_q[n]_q[n]_{q^2}. \tag{25}$$

For

$$\sum_{j=0}^{n-1}[n+1+2j]_q\, q^{n-1-j}=\frac{1}{1-q}\left(\sum_{j=0}^{n-1}q^{n-1-j}-\sum_{j=0}^{n-1}q^{2n+j}\right)=\frac{(1-q^n)(1-q^{2n})}{(1-q)^2}=(1+q)[n]_q[n]_{q^2}.$$

### 4.2

A similar identity is

$$\sum_{k=0}^{n}[mn+k]_{q^2}\, q^{n-k}=\begin{bmatrix}n+1\\2\end{bmatrix}_q[2m+1]_{q^n}. \tag{26}$$

The left-hand side is

$$\sum_{k=0}^{n}q^{n-k}\frac{1-q^{2mn+2k}}{1-q^2}=\frac{1}{1-q^2}\left(\sum_{k=0}^{n}q^{n-k}-q^{(2m+1)n+k}\right)=\frac{[n+1]_q}{1+q}\frac{1-q^{(2m+1)n}}{1-q}=\frac{[n+1]_q[n]_q}{1+q}[2m+1]_{q^n}.$$

### 4.3

C. Wheatstone observed that

$$\sum_{j=\frac{3^n+1}{2}}^{\frac{3^{n+1}-1}{2}}j=3^{2n} \tag{27}$$

implies

$$\sum_{j=1}^{\frac{3^{n+1}-1}{2}}j=\sum_{k=0}^{n}3^{2k}. \tag{28}$$

A $q-$ analogue is

**Theorem 3**

$$\sum_{j=\frac{3^n+1}{2}}^{\frac{3^{n+1}-1}{2}} [j]_{q^2} q^{\frac{3^{n+1}-1-2j}{2}} = \left[3^n\right]_q \left[3^n\right]_{q^2} \tag{29}$$

*and therefore*

$$\sum_{j=1}^{\frac{3^{n+1}-1}{2}} q^{\frac{3^{n+1}-2j-1}{2}} [j]_{q^2} = \sum_{k=0}^{n} q^{\frac{3^{n+1}-3^{k+1}}{2}} \left[3^k\right]_q \left[3^k\right]_{q^2} = \begin{bmatrix} \frac{3^{n+1}+1}{2} \\ 2 \end{bmatrix}_q. \tag{30}$$

**Proof**

Observe that (cf. [3])

$$\sum_{k=0}^{n} [k]_{q^2} q^{n-k} = \begin{bmatrix} n+1 \\ 2 \end{bmatrix}_q. \tag{31}$$

This is the special case $m=0$ of (26).

For $a(n) = \dfrac{3^n - 1}{2}$ the left-hand side of (29) becomes

$$\begin{bmatrix} a(n+1)+1 \\ 2 \end{bmatrix}_q - q^{a(n+1)-a(n)} \begin{bmatrix} a(n)+1 \\ 2 \end{bmatrix}_q = \frac{\left(1-q^{a(n+1)+1}\right)\left(1-q^{a(n+1)}\right) - q^{a(n+1)-a(n)}\left(1-q^{a(n)+1}\right)\left(1-q^{a(n)}\right)}{(1-q)(1-q^2)}$$

$$= \frac{\left(1-q^{a(n+1)+1}\right)\left(1-q^{a(n+1)}\right) - q^{a(n+1)-a(n)}\left(1-q^{a(n)+1}\right)\left(1-q^{a(n)}\right)}{(1-q)(1-q^2)} = \frac{\left(1-q^{a(n+1)-a(n)}\right)\left(1-q^{a(n+1)+a(n)+1}\right)}{(1-q)(1-q^2)}.$$

This gives

$$\sum_{j=\frac{3^n+1}{2}}^{\frac{3^{n+1}-1}{2}} [j]_{q^2} q^{\frac{3^{n+1}-1-2j}{2}} = \frac{\left(1-q^{3^n}\right)\left(1-q^{2 \cdot 3^n}\right)}{(1-q)(1-q^2)} = \left[3^n\right]_q \left[3^n\right]_{q^2} = \frac{\left(1-q^{3^n}\right)^2 \left(1+q^{3^n}\right)}{(1-q)(1-q^2)} = \left[3^n\right]_q^2 \frac{1+q^{3^n}}{1+q}.$$

and therefore $\displaystyle\sum_{j=1}^{\frac{3^{n+1}-1}{2}} q^{\frac{3^{n+1}-2j-1}{2}} [j]_{q^2} = \sum_{k=0}^{n} q^{\frac{3^{n+1}-3^{k+1}}{2}} \left[3^k\right]_q \left[3^k\right]_{q^2}$, which by (31) implies (30).

**4.4**

A $q-$ analogue of $\sum_{j=0}^{n-1}(2(n+1)j+1) = n^3$ is

$$\sum_{j=0}^{n-1}[2(n+1)j+1]_q q^{n^2-1-(n+1)j} = [n^2]_q [n]_{q^{n+1}}. \tag{32}$$

A $q-$ analogue of $\sum_{j=0}^{n-1}(2k+1)n = n^3$ is

$$\sum_{k=0}^{n-1}[(2k+1)n]_q q^{n(n-1)-nk} = [n^2]_q [n]_{q^n}. \tag{33}$$

More generally we get for $m \in \mathbb{N}$

$$\sum_{k=0}^{n-1}[(2k+1)n^m]_q q^{n^m(n-1)-kn^m} = [n^{m+1}]_q [n]_{q^{n^m}}. \tag{34}$$

The proofs are straightforward and will be omitted.

**4.5**

C. Wheatstone also observed that

$$1^3 = 1$$
$$2^3 = 3+5$$
$$3^3 = 6+9+12$$
$$4^3 = 10+14+18+22$$
$$5^3 = 15+20+25+30+35$$
$$\ldots\ldots\ldots$$

More generally this gives $\sum_{j=0}^{n-1}\left(\binom{n+1}{2}+jn\right) = n^3$.

We will now prove three different $q-$ analogues of this identity.

**Theorem 4**

$$\sum_{j=0}^{n-1}\left[\binom{n+1}{2}+nj\right]_q q^{\binom{n}{2}-\binom{j+1}{2}} = [n^2]_q \sum_{j=0}^{n-1} q^{\binom{n}{2}-\binom{j+1}{2}}, \tag{35}$$

$$\sum_{j=0}^{n-1}\left[\binom{n+1}{2}+nj\right]_{q^2} q^{n(n-1)-nj} = [n^2]_{q^2} [n]_{q^n}, \tag{36}$$

*and*

$$\sum_{j=0}^{n-1}\left(q^{n-1-j}\left[\binom{n+1}{2}\right]_q + q^{2n-j}[n]_q [j]_{q^2}\right) = [n]_q^2 [n]_{q^2}. \tag{37}$$

**Proof**

(35) follows from

$$\frac{1}{1-q}\sum_{j=0}^{n-1}\left(q^{\binom{n}{2}-\binom{j+1}{2}}\left(1-q^{\binom{n+1}{2}+nj}\right)\right) = \frac{q^{\binom{n}{2}}}{1-q}\left(\sum_{j=0}^{n-1} q^{-\binom{j+1}{2}} - q^{\binom{n+1}{2}}\sum_{j=0}^{n-1} q^{nj-\binom{j+1}{2}}\right)$$

$$= \frac{q^{\binom{n}{2}}}{1-q}\left(\sum_{j=0}^{n-1} q^{-\binom{j+1}{2}} - q^{\binom{n+1}{2}}\sum_{j=0}^{n-1} q^{\binom{n}{2}-\binom{n-j}{2}}\right) = \frac{1}{1-q}\sum_{j=0}^{n-1} q^{\binom{n}{2}-\binom{j+1}{2}}\left(1-q^{n^2}\right).$$

For the second identity the left-hand side is

$$\frac{1}{1-q^2}\sum_{j=0}^{n-1}\left(q^{n^2-n-nj}\left(1-q^{2\binom{n+1}{2}+2nj}\right)\right) = \frac{1}{1-q^2}\left(\sum_{j=0}^{n-1} q^{nj} - q^{2n^2}\sum_{j=0}^{n-1} q^{nj}\right) = \frac{(1-q^{n^2})(1-q^{2n^2})}{(1-q^n)(1-q^2)},$$

which implies the right-hand side.

After dividing (37) by $[n]_q$ and multiplying with $(1-q)(1-q^2)$ the left-hand side is

$$(1-q)\sum_{j=0}^{n-1}\left(q^{n-1-j}(1-q^{n+1}) + q^{2n-j}(1-q^{2j})\right) = (1-q)\sum_{j=0}^{n-1}\left(q^{n-1-j} - q^{2n-j} + q^{2n-j} - q^{2n+j}\right)$$

$$= (1-q)\left(\sum_{j=0}^{n-1} q^j - q^{n+1}\sum_{j=0}^{n-1} q^j + (q^{n+1}-q^{2n})\sum_{j=0}^{n-1} q^j\right) = (1-q^n)(1-q^{2n}).$$

This proves (37).

## 4.6

Finally we look for a $q-$ analogue of $\sum_{j=0}^{2n}\left((2n-1)^2+8j\right)=(2n+1)^3$.

**Theorem 5**

$$\sum_{j=0}^{2n}\left[(2n-1)^2+8j\right]_q q^{8n-4j} = [2n+1]_q [2n+1]_{q^4} [2n+1]_{q^{2n+1}}. \tag{38}$$

**Proof**

$$\frac{1}{1-q}\sum_{j=0}^{2n} q^{8n-4j}\left(1-q^{(2n-1)^2+8j}\right) = \frac{1}{1-q}\sum_{j=0}^{2n} q^{4j} - \frac{q^{(2n-1)^2+8n}}{1-q}\sum_{j=0}^{2n} q^{4j} = \frac{\left(1-q^{4(2n+1)}\right)\left(1-q^{(2n+1)^2}\right)}{(1-q)(1-q^4)}$$

$$= [2n+1]_q [2n+1]_{q^4} [2n+1]_{q^{2n+1}}.$$